\documentclass[twoside,12pt]{amsart}
\usepackage{fullpage}
\usepackage{amscd}
\newcounter{list_count}
\newenvironment{rlist}{\begin{list}
		{(\roman{list_count})}{\usecounter{list_count}
		\setlength{\rightmargin}{\leftmargin}}}{\end{list}}

\begin{document} 
%\begin{spacing}{1}
%\input{psbox.tex}
%\psfordvips
%%%%%%%%%%%%%%%%%%%%%%%%%%%%%%%%%%%%%%%%%%%%%%%%%%%%%%%%%%%%%%%%%%%%%%
\setcounter{tocdepth}{2}
%\def\contentsname{}
%%%%%%%%%%%%%%%%%%%%%%%%%%%%%%%%%%%%%%%%%%%%%%%%%%%%%%%%%%%%%%%%%%%%%%
\makeatletter
\def\@biblabel#1{#1.}%
\newenvironment{thelocalbibliography}[1]{%
%  \@xp\section\@xp*\@xp{\refname}%
%  \normalfont\footnotesize
\normalfont\small\labelsep .5em\relax
  \renewcommand\theenumiv{\arabic{enumiv}}\let\p@enumiv\@empty
  \list{\@biblabel{\theenumiv}}{\settowidth\labelwidth{\@biblabel{#1}}%
    \leftmargin\labelwidth \advance\leftmargin\labelsep
    \usecounter{enumiv}}%
  \sloppy \clubpenalty\@M \widowpenalty\clubpenalty
  \sfcode`\.=\@m
}{%
  \def\@noitemerr{\@latex@warning{Empty `thebibliography' environment}}%
  \endlist
}
\makeatother
%%%%%%%%%%%%%%%%%%%%%%%%%%%%%%%%%%%%%%%%%%%%%%%%%%%%%%%%%%%%%%%%%%%%%%
\outer\def\refinsection#1{%{
\unskip\nobreak\hfil\penalty50
         \hskip2em\hbox{}\nobreak\hfil(\S{\it\ref{#1}\/})
         \parfillskip0pt\finalhyphendemerits=0\par}% TeXbook, p 106.
%    {\ifdim\lastskip<\bigskipamount\removelastskip\penalty-10\bigskip\fi}}
%%%%%%%%%%%%%%%%%%%%%%%%%%%%%%%%%%%%%%%%%%%%%%%%%%%%%%%%%%%%%%%%%%%%%%
\newcommand{\FAPESP}{FAPESP}
\newcommand{\CNPq}{CNPq}
%%%%%%%%%%%%%%%%%%%%%%%%%%%%%%%%%%%%%%%%%%%%%%%%%%%%%%%%%%%%%%%%%%%%%%
\def\({\left(}
\def\){\right)}
\def\<{\left\langle}
\def\>{\right\rangle}
\def\:{\colon}
%%%%%%%%%%%%%%%%%%%%%%%%%%%%%%%%%%%%%%%%%%%%%%%%%%%%%%%%%%%%%%%%%%%%%%
\def\choose{\atopwithdelims()}
%%%%%%%%%%%%%%%%%%%%%%%%%%%%%%%%%%%%%%%%%%%%%%%%%%%%%%%%%%%%%%%%%%%%%%
\def\NP{\mathcal{NP}}
\def\P{\mathcal{P}}
\def\RR{\mathcal{R}}
\def\complexes{\mathbb{C}}
\def\reals{\mathbb{R}}
%%%%%%%%%%%%%%%%%%%%%%%%%%%%%%%%%%%%%%%%%%%%%%%%%%%%%%%%%%%%%%%%%%%%%%
\def\arrow{\to}
\def\indarrow{{%\displaystyle
  \mathrel{\mathop{\longrightarrow}\limits^{\hbox{\tiny ind}}}}}
\def\orarrow{{%\displaystyle
  \mathrel{\mathop{\longrightarrow}\limits^{\hbox{\tiny or}}}}}
\def\aorarrow{{%\displaystyle
 of all real continuous  functions on $I$

  \mathrel{\mathop{\longrightarrow}\limits^{\hbox{\tiny acy}}}}}
\def\mcarrow{{%\displaystyle
  \mathrel{\mathop{\longrightarrow}\limits^{\hbox{\tiny mc}}}}}
\def\rind{r^{\rm ind}}
\def\re{r_{\rm e}}
\def\rinde{r_{\rm e}^{\rm ind}}
\def\orr{r^{\rm or}}
\def\aorr{r_{\rm a}^{\rm or}}
\def\ghat{\hat g}
\def\ex{\mathop{\rm ex}\nolimits}
\def\find{f_{\rm ind}}
\def\gind{g_{\rm ind}}
\def\cP{\mathcal{P}}
%%%%%%%%%%%%%%%%%%%%%%%%%%%%%%%%%%%%%%%%%%%%%%%%%%%%%%%%%%%%%%%%%%%%%%
%\setlength{\textwidth}{5cm}
\theoremstyle{plain}
\newtheorem{theorem}{Theorem}%[section]
\newtheorem{problem}[theorem]{Problem}
\newtheorem{definition}{Definition}
\newtheorem{conjecture}[theorem]{Conjecture}
\newtheorem{corollary}[theorem]{Corollary}
\newtheorem{lemma}[theorem]{Lemma}
\newtheorem{proposition}[theorem]{Proposition}
%%%%%%%%%%%%%%%%%%%%%%%%%%%%%%%%%%%%%%%%%%%%%%%%%%%%%%%%%%%%%%%%%%%%%%
%\makeatletter
%\newcommand{\subsubsubsection}[1]{\smallskip\noindent\textit{#1}\enspace} 
%\makeatother
%%%%%%%%%%%%%%%%%%%%%%%%%%%%%%%%%%%%%%%%%%%%%%%%%%%%%%%%%%%%%%%%%%%%%%
\def\IMSmarkvadjust{-8 pt}
\def\IMSmarkhadjust{0 pt}
\def\SBIMSMark#1#2#3{
 \font\SBF=cmss10 at 10 true pt
 \font\SBI=cmssi10 at 10 true pt
 \setbox0=\hbox{\SBF Stony Brook IMS Preprint \##1}
 \setbox2=\hbox to \wd0{\hfil \SBI #2}
 \setbox4=\hbox to \wd0{\hfil \SBI #3}
 \setbox6=\hbox to \wd0{\hss
             \vbox{\hsize=\wd0 \parskip=0pt \baselineskip=10 true pt
                   \copy0 \break%
                   \copy2 \break% 
                   \copy4 \break}}
 \dimen0=\ht6   \advance\dimen0 by \vsize \advance\dimen0 by 8 true pt
                \advance\dimen0 by -\pagetotal
	        \advance\dimen0 by \IMSmarkvadjust
 \dimen2=\hsize \advance\dimen2 by .25 true in
	        \advance\dimen2 by \IMSmarkhadjust

%
%   Check for publication info
%
%  \newread\jref
  \openin2=publishd.tex
  \ifeof2\setbox0=\hbox to 0pt{}
  \else 
     \setbox0=\hbox to 3.1 true in{
                \vbox to \ht6{\hsize=3 true in \parskip=0pt  \noindent  
                {\SBI Published in modified form:}\hfil\break
                \input publishd.tex 
                \vfill}}
  \fi
  \closein2
  \ht0=0pt \dp0=0pt
 \ht6=0pt \dp6=0pt
 \setbox8=\vbox to \dimen0{\vfill \hbox to \dimen2{\copy0 \hss \copy6}}
 \ht8=0pt \dp8=0pt \wd8=0pt
 \copy8
 \message{*** Stony Brook IMS Preprint #1, #2. #3 ***}
}
%%%%%%%%%%%%%%%%%%%%%%%%%%%%%%%%%%%%%%%%%%%%%%%%%%%%%%%%%%%%%%%%%%%
\renewcommand{\today}{\relax}
\pagestyle{myheadings}
\markboth{\sc \hfill [Rigidity of $C^2$ infinitely renormalizable unimodal
                      maps \hfill \quad} 
         {\sc \quad\hfill W. de Melo and A. Pinto \hfill}
   
\title[Rigidity of
$C^2$ infinitely renormalizable unimodal maps]
{Rigidity of
$C^2$ infinitely renormalizable unimodal maps}       

\bigskip

\author{W. de Melo}
\address[W. de Melo]{IMPA,
Rio de Janeiro, Brazil.}
%% Note the doubled @@:
\email[W. de Melo]{demelo@impa.br}

\author{A. A. Pinto}
\address[A. A. Pinto]{DMA, Faculdade de Ci\^encias, Universidade
do Porto\\ 4000 Porto, Portugal.}
%% Note the doubled @:
\email[A. A. Pinto]{aapinto@fc.up.pt}

%\address{}
%\email{}
%\thanks{This work has been supported by PRONEX-Dynamical Systems}

\begin{abstract}
Given     $C^2$
infinitely renormalizable unimodal maps 
$f$ and $g$ with a quadratic critical point and 
the same bounded combinatorial type,
we prove that
they are  $C^{1+\alpha}$  conjugate 
along the closure of the corresponding 
forward  orbits of the critical points,
for some $\alpha>0$.
\end{abstract}


\maketitle
\tableofcontents
\thispagestyle{empty}

\SBIMSMark{1999/6}{May 1999}{}

\newpage

\section{Introduction}
 \label{jhvbmnk}

It was already clear more than 20 years ago, from the work of Coullet-Tresser 
and Feigenbaum, that the small scale geometric properties of the orbits of 
some one dimensional dynamical systems were related to the dynamical 
behavior of a non-linear operator, the renormalization operator, acting on 
a space of dynamical systems. This conjectural picture was mathematically 
established for some classes of analytic maps by Sullivan, McMullen and 
Lyubich. Here we will extend this description to the space of $C^2$ maps and 
prove a rigidity result for
a class of unimodal maps of the interval. 
As it is well-known, a unimodal map is a smooth endomorphism of a 
compact interval that has a unique critical point
which is a turning point. Such a map is renormalizable if there exists 
an interval neighborhood of the critical point such that
the first return map to this interval is again a unimodal map, and 
the return time is greater than one. The map
is infinitely renormalizable if there exist  
 such intervals with arbitrarily high return times.
We say that two maps have the same combinatorial type if the map that sends the 
i-th iterate of the critical point of the first map into 
the i-th iterate of the critical point of the second map,
for all $i \ge 0$,
 is order preserving.
Finally, we say that the combinatorial type of an 
infinitely renormalizable map is bounded if 
the ratio of any two  consecutive return times is uniformly bounded.

A unimodal map $f$ is $C^r$ with a quadratic critical point 
if 
$f= \phi_f \circ p \circ \psi_f$, where  $p(x)=x^2$ and $\phi_f$, $\psi_f$ are $C^r$
diffeomorphisms. 
Let
 $c_f$ be the critical point of $f$.
 In this paper we will prove the following rigidity result.

\begin{theorem}
\label{aaa}
Let $f$ and
 $g$ be  $C^2$ unimodal maps with a quadratic 
critical point which are infinitely renormalizable and have the same
bounded combinatorial type. Then  there exists a $C^{1+\alpha}$
diffeomorphism $h$ of the real line such that 
$h(f^i(c_f))=g^i(h(c_g))$ for every integer $i\ge 0$.
\end{theorem}

 We observe that in Theorem \ref{aaa}
 the H\"older exponent $\alpha>0$ depends only upon the bound 
of the combinatorial type of the maps $f$ and $g$.
Furthermore,as we will see in Section \ref{sec2}, the maps $f$ and $g$ are
 smoothly conjugated to $C^2$ normalized unimodal maps 
 $F=\phi_F \circ p$ and $G=\phi_G \circ p$
with  critical  value $1$,
and the  H\"older constant for the smooth conjugacy 
between the normalized maps $F$ and $G$ depends only upon 
the combinatorial type of $F$ and $G$, and
 upon the norms $||\phi_F||_{C^2}$ and $||\phi_G||_{C^2}$.

The conclusion of the above rigidity theorem was first obtained by
McMullen in \cite{mullen2}  
 under the extra hypothesis that  $f$ and $g$
extend to quadratic-like maps in neighborhoods of the dynamical intervals
in the complex plane. Combining this last statement  with the complex bounds 
of Levin and van Strien in \cite{levin}, we get the 
existence of a $C^{1+\alpha}$ map $h$ which is a conjugacy
along 
the critical orbits for  infinitely renormalizable real analytic maps
with the same
bounded combinatorial type. 
 We   extended this result to $C^2$ unimodal maps in 
 Theorem  \ref{aaa},  by combining many results and ideas of Sullivan
in \cite{S1} 
with recent results of McMullen in \cite{mullen1}, 
  in \cite{mullen2}, and  of Lyubich in \cite{lyubich}
on the hyperbolicity of 
the renormalization operator $R$
(see the definition of $R$ in the next section). 
A main lemma used in the proof of Theorem \ref{aaa}
is the following:

\begin{lemma}
\label{xzx} 
Let $f$ be a $C^2$ infinitely renormalizable 
map with bounded combinatorial type.
Then  there exist positive constants 
$\eta<1$,  $\mu$ and $C$, 
and a real quadratic-like map  $f_n$ 
with conformal modulus  greater than or equal to $\mu$,
and with the same combinatorial type as the $n$-th
renormalization $R^nf$ of $f$ such that 
$$||R^nf-f_n||_{C^0} < C \eta^n$$
for every $n \ge 0$.
\end{lemma} 

We observe that in this lemma,
the positive constants $\eta<1$ and  $\mu$
  depend only upon    the bound 
of the combinatorial type of the map  $f$.
For normalized unimodal maps $f$,
the positive constant $C$ 
  depends only upon    the bound 
of the combinatorial type of the map  $f$
and upon the norm $||\phi_f||_{C^2}$. 

This lemma 
 generalizes a Theorem of Sullivan
 (transcribed as Theorem \ref{th3} in Section \ref{sec2})
by adding that the map $f_n$ has the  same combinatorial type as the $n$-th
renormalization $R^nf$ of $f$.

Now, let us   describe the proof of Theorem \ref{aaa}
which also  shows the relevance of Lemma \ref{xzx}:
 let $f$ and $g$
be  $C^2$   
infinitely renormalizable unimodal maps
with the same bounded combinatorial type.
Take $m$ to be of the order of a large but fixed fraction of $n$,
and note that $n-m$ is also a fixed fraction of $n$.
By Lemma \ref{xzx}, we obtain a real
quadratic-like map  $f_m$
exponentially close to $R^mf$, 
and a real 
quadratic-like map  $g_m$
exponentially close to $R^mg$.
Then  we use Lemma \ref{l1} of Section \ref{leee} to prove 
 that the renormalization $(n-m)$-th iterates 
$R^nf$ of $R^mf$, and $R^ng$ of $R^mg$
 stay exponentially close to the 
 $(n-m)$-th iterates   
$R^{n-m}f_m$ of $f_m$  and 
$R^{n-m}g_m$ of $g_m$, respectively.
Again, by Lemma \ref{xzx}, we have 
that $f_m$ and $g_m$ have conformal modulus 
universally bounded away from zero,
and     have the same bounded combinatorial type of $R^mf$ and $R^mg$. 
Thus, by the main result  of 
McMullen in \cite{mullen2},
 the renormalization $(n-m)$-th iterates $R^{n-m}f_m$ of
   $f_m$ and $R^{n-m}g_m$ of $g_m$ are
     exponentially close.
 Therefore, 
$R^nf$ is exponentially close to $R^{n-m}f_m$, 
$R^{n-m}f_m$ is exponentially close to $R^{n-m}g_m$,
and $R^{n-m}g_m$  is exponentially close to $R^ng$, 
and so, by the triangle inequality, 
the $n$-th iterates  $R^nf$ of $f$ and $R^ng$ of $g$
  converge  exponentially fast to each other. Finally, 
by   Theorem 9.4 in the book \cite{vSdM}  of de Melo and van Strien,  
 we conclude  that
    $f$ and $g$ are $C^{1+\alpha}$ 
  conjugate  along the closure of their critical orbits.

Let us point out the main
ideas in the proof of Lemma \ref{xzx}:
Sullivan   in \cite{S1}   proves that 
 $R^nf$   is exponentially close to a
 quadratic-like map  $F_n$ which  has 
conformal modulus universally bounded away from zero. 
The quadratic-like map  $F_n$ 
determines   a unique quadratic  map  $P_{c(F_n)}(z)=1-c(F_n)z^2$ which is
hybrid
conjugated to $F_n$ by a $K$  quasiconformal homeomorphism, 
where $K$ depends only upon the conformal modulus of $F_n$ 
(see Theorem  1 of   Douady-Hubbard  in \cite{douady},   and Lemma \ref{auxfinding} in Section \ref{llfff}).
In \cite{lyubich}, Lyubich proves the bounded geometry 
of the Cantor set consisting of all the parameters 
of the quadratic family $P_c(z)=1-c z^2$
corresponding to infinitely renormalizable maps 
with  combinatorial type bounded by $N$
 (see definition in Section \ref{sec2} and the proof of Lemma \ref{xzx}).
In Lemma \ref{c3} of Section \ref{leee}, we prove that $R^nf$ and $F_n$  
have exponentially close renormalization types. 
Therefore, letting $c_n$ be the parameter corresponding 
to the quadratic map $P_{c_n}$ with 
the same combinatorial type as $R^nf$,
we have, from the above result of Lyubich,
  that $c(F_n)$ and $c_n$ are exponentially close.
In Lemma \ref{finding} of Section \ref{llfff}, we use   holomorphic motions
to prove the existence of a real quadratic-like map   $f_n$ 
which is hybrid conjugated to $P_{c_n}$, and 
has the following essential property: 
the distance between $F_n$ and $f_n$ is proportional to the 
distance between $c(F_n)$ and $c_n$ raised to some positive constant. 
Therefore, the  real quadratic-like map 
$f_n$ has the same combinatorial type as $R^nf$, and 
  $f_n$ is exponentially close to $F_n$.
Since the map  $F_n$ is   exponentially close to $R^nf$,
we obtain that the map  $f_n$ is also exponentially close to $R^nf$.

The example of Faria and de Melo in \cite{welingtonfaria} 
for critical circle maps can be adapted to prove the
existence of a pair of $C^\infty$ unimodal maps, with the same unbounded
combinatorial 
type,  such that the conjugacy $h$  has no $C^{1+\alpha}$ 
extension to the reals for any $\alpha>0$.

\section{Shadowing unimodal maps}
\label{sec2}

 A  $C^r$  unimodal map
 $F:I \to I$ is  {\it normalized}
if  $I=[-1,1]$,  $F=\phi_F \circ p$,
 $F(0) = 1$,  and  $\phi_F:[0,1] \to I$ is  a $C^r$ diffeomorphism.
A $C^r$ unimodal map  $f= \phi_f \circ p \circ \psi_f$ 
with  quadratic critical point  either has trivial dynamics 
or has an invariant interval where it is $C^r$ 
conjugated to a  $C^r$ normalized unimodal map
 $F$. Take, for instance, the map
$$
\phi_F(x) =
\left( \psi_f^{-1} \circ \phi_f (0) \right)^{-1} \cdot
 \psi_f^{-1} \circ \phi_f 
\left(
\left( \psi_f^{-1} \circ \phi_f (0) \right)^{-2} \cdot x \right) \ .
$$
 Therefore, from now on we will only
consider $C^r$ normalized unimodal maps $f$.

The map $f$ is {\it renormalizable} if there is
a closed interval $J$ centered at the origin,
strictly contained in $I$, and $l > 1$ such that
the intervals $J, \ldots, f^{l-1}(J)$
are disjoint, $f^l(J)$ is strictly contained in $J$ and 
$f^l(0) \in \partial J$.
If $f$ is renormalizable, we always consider the smallest $l>1$  and
the minimal interval  $J_f=J$  with the above
properties. The set of all renormalizable maps is 
an open set in the $C^0$ topology.
 The {\it renormalization operator $R$} acts on
renormalizable maps  $f$ by
 $Rf=\psi \circ f^l \circ \psi^{-1}:I \to I$,
where  $\psi:J_f \to I$ is the restriction
of a linear map sending $f^l(0)$ into $1$.
Inductively, the map $f$ is {\it $n$ times renormalizable} 
  if  $R^{n-1}f$  is renormalizable. 
 If $f$ is $n$ times renormalizable for
every $n > 0$, then $f$ is    {\it infinitely renormalizable}.

Let $f$ be a renormalizable map.
We label   the intervals  $J_f, \ldots, f^{l-1}(J_f)$
of $f$
by $1,\ldots,l$
 according to their embedding on
the real line, from the left to the right.
{\it The  permutation}
$\sigma_f:\{1,\ldots,l\} \to \{1,\ldots,l\}$
is defined by
$\sigma_f(i)=j$ if the interval labeled by $i$ is mapped by
$f$ to the interval labeled by $j$.
The {\it renormalization type} of
 an $n$ times renormalizable map  $f$
is given by the sequence $\sigma_f, \ldots,\sigma_{R^n f}$.
An $n$ times renormalizable map  $f$  has
{\it renormalization type bounded by $N > 1$} 
if the number of elements of the domain of each permutation $\sigma_{R^mf}$
is less than or equal to $N$ for every $0 \le m \le n$.
We have the analogous notions for  infinitely renormalizable maps.

Note that if any two maps are $n$ times renormalizable 
and have the same combinatorial type 
(see definition in the introduction),
then they have the  same renormalization type.
The converse is also true in the case of infinitely 
renormalizable maps.  An infinitely 
renormalizable map   has combinatorial  type bounded by $N > 1$
if the  renormalization   type is   bounded by $N$.

If $f=\phi \circ p$ is $n$ times renormalizable, 
and $\phi \in C^2$,  
there is a $C^2$ diffeomorphism $\phi_n$
satisfying $R^nf=\phi_n \circ p$.
The {\it nonlinearity} ${\rm nl}(\phi_n)$ of  $\phi_n$ is defined by
 $${\rm nl}(\phi_n)=\sup_{x \in p(I)}
\left| \frac{\phi_n''(x)}{\phi_n'(x)} \right| \ .
$$

Let $\mathcal{I}(N,b)$ be the set of all
  $C^2$ normalized unimodal maps $f=\phi \circ p$
with the following properties:
\begin{rlist}
\item $f$ is
infinitely renormalizable;
\item  the combinatorial type of $f$ is
 bounded by $N$;
\item $||\phi||_{C^2} \le b$.
\end{rlist}
 
\begin{theorem}
\label{t2}
{\bf (Sullivan \cite{S1})}
There exist positive  constants  $B$ and $n_1(b)$ such that, 
for every $f \in \mathcal{I}(N,b)$,    the $n$-th renormalization
 $R^nf=\phi_n \circ p$ of $f$ has the property that
${\rm nl}(\phi_n) \le B$  for every
$n \ge n_1$.
\end{theorem} 
 
This theorem together with Arzel\'a-Ascoli's Theorem 
implies that, for every   $0 \le \beta < 2$, and for every $n  \ge n_1(b)$,
 the  renormalization iterates  $R^nf$
 are contained in a compact set
of unimodal maps with respect to the
 $C^{\beta}$ norm. We will use this 
fact in the  proof of Lemma \ref{l2} below.

\subsection{Quadratic-like maps}

A {\it quadratic-like map $f:V \to W$}
 is a holomorphic map with the property that $V$ and $W$
are simply connected domains with the closure
of $V$ contained in   $W$, and $f$ is a
 degree two branched covering map. We add an extra condition that $f$ has
 a continuous extension to the boundary of $V$.
The {\it conformal modulus of a quadratic-like map} $f:V \to W$
is equal to the conformal modulus of the annulus
$W \setminus \overline{V}$. A {\it real quadratic-like
map} is a quadratic-like map which
commutes with  complex conjugation.

The {\it  filled  Julia set ${\mathcal K}(f)$ of $f$}
is the set
$\{z:f^n(z) \in V, ~{\rm for} ~{\rm all}~n \ge 0\}$.
Its boundary is the {\it Julia set} ${\mathcal J}(f)$ of $f$.
These sets ${\mathcal J}(f)$ and ${\mathcal K}(f)$ are connected if 
the critical point of $f$ is contained in ${\mathcal K}(f)$.

Let ${\mathcal Q}(\mu)$ be the set of all  real quadratic-like
maps  $f:V \to W$   satisfying the following properties:
\begin{rlist}
\item the Julia set  ${\mathcal J}(f)$ of $f$ is connected;
\item the  conformal modulus of $f$ is
  greater than or equal to  $\mu$, and less than or equal to
$2\mu$;
\item $f$ is normalized to have the critical point at the origin, and the critical value
at one.
\end{rlist}
By Theorem 5.8 in page 72 of \cite{mullen1}, 
the set ${\mathcal Q}(\mu)$ is compact in the  Carath\'eodory topology
taking the critical point as the base point
(see definition in page 67 of  \cite{mullen1}).

\begin{theorem}
\label{th3} 
{\bf (Sullivan \cite{S1})}
There exist positive constants $\gamma(N) < 1$, 
 $C(b,N)$, and $\mu(N)$   
 with the following property:
if $f \in \mathcal{I}(N,b)$, then   
there exists   $f_n \in {\mathcal Q}(\mu)$
such that
$||R^nf-f_n||_{C^0} \le C \gamma^n$. 
\end{theorem}

In the following sections,
 we will develop the results
that  will be used in the last section to prove 
 the generalization of Theorem \ref{th3} 
(as stated in Lemma \ref{xzx}), 
 and  to prove Theorem \ref{aaa}.

\subsection{Maps with   close combinatorics} 
 \label{leee}

Let $D(\sigma)$ be the open set of all
$C^0$
renormalizable unimodal maps $f$
 with renormalization  type $\sigma_f=\sigma$. 
The open 
sets $D(\sigma)$ are pairwise disjoint.
Let $E(\sigma)$ be the complement  of $D(\sigma)$
in the set of all $C^0$ unimodal maps $f$.

\begin{lemma}
\label{l2}
There exist positive constants $n_2(b)$ and $\epsilon(N)$
with the following property:
for every  $f \in \mathcal{I}(N,b)$,
for every 
 $n \ge n_2$, and for every 
   $g \in  E(\sigma_{R^nf})$,  we have
 $
||R^nf- g||_{C^0} > \epsilon .
$
\end{lemma}

 \begin{proof}
Suppose, by contradiction, that
there is a sequence
$R^{m_1}f_1$, $R^{m_2}f_2$, $\ldots$
with the property that for a chosen $\sigma$
there is a sequence $g_1,g_2, \ldots \in  E(\sigma)$
satisfying
 $||R^{m_i}f_i- g_i||_{C^0} < 1/i$. 
By Theorem \ref{t2}, there are  $B>0$ and $n_1(b)\ge 1$ such that
the maps $R^{m_i}f_i$
have  nonlinearity bounded by $B>0$ for all $m_i \ge n_1$.
By Arzela-Ascoli's Theorem, there is a   subsequence
$R^{m_{i_1}}f_{i_1},R^{m_{i_2}}f_{i_2},\ldots$ which converges in the
$C^0$
topology
  to a   map $g$. Hence, the map $g$ is contained
 in the boundary of $E(\sigma)$
and is infinitely renormalizable.
However, a map contained
 in the boundary of $E(\sigma)$
  is not renormalizable, and so we get   a contradiction.
\end{proof}

\begin{lemma}
\label{l1} 
There exist positive constants $n_3(N,b)$ and $L(N)$
with the following property:
for every  $f \in \mathcal{I}(N,b)$, 
for every $C^2$ renormalizable unimodal map 
$g$, and for every  $n>n_3$, 
we have
$$||R^nf-Rg||_{C^0} \le L  ||R^{n-1}f-g||_{C^0} \ .$$
\end{lemma}

\begin{proof}
In the proof of this lemma we will use
the inequality (\ref{qwe1}) below.
 Let $f_1,\ldots,f_m$ 
be maps with $C^1$ norm bounded by some constant
$d > 0$, and let $g_1,\ldots,g_m$ be $C^0$ maps.
By induction on $m$, and
by the Mean Value Theorem, 
there is $c(m,d) >0$ such that
\begin{equation}
\label{qwe1}
||f_1 \circ \ldots \circ f_m
- g_1 \circ \ldots \circ g_m ||_{C^0}
\le c \max_{i=1,\ldots,m} \{||f_i - g_i||_{C^0}\} \ .
\end{equation}

Set $n_3=\max \{n_1,n_2\}$, where $n_1(b)$ is defined as 
in Theorem \ref{t2}, and  $n_2(b)$ is defined as 
in Lemma \ref{l2}.
Set $F=R^{n-1}f$ with $n \ge n_3$.
We start by considering  the simple case (a),  where 
 $F$ and $g$ do not have the same 
renormalization type, and conclude with the complementary case (b).
In case (a), by Lemma \ref{l2},
there is  $\epsilon(N)>0$ with the property that
$$
||RF-Rg||_{C^0}
  \le 2 \le 2 \epsilon^{-1} ||F-g||_{C^0} \ .$$
In case (b),  there is $1 < m \le N$ such that
 $R F(x)= a_F F^m (a_F^{-1} x)$, 
and  $R g(x)= a_g  g^m (a_g^{-1}x)$,
where  $a_F= F^m(0)$ and $a_g=g^m(0)$.
By Theorem \ref{t2}, 
there is a positive constant $B(N)$
bounding the nonlinearity of $F$. 
Since  the set of all infinitely renormalizable unimodal maps  $F$
with nonlinearity bounded by $B$ is a compact set
with respect to the $C^0$ topology,
and since $a_F$  
varies continuously  with $F$, 
 there is $S(N)>0$ with the property that  $|a_F| \ge S$.
Again, by Theorem \ref{t2}, and by inequality (\ref{qwe1}),
there is $c_1(N) > 0$   such that
\begin{equation}
\label{qwe2}
||F^m-g^m||_{C^0} 
\le c_1 ||F-g||_{C^0} \ .
\end{equation}
Thus, 
\begin{equation}
\label{qwe233}
|a_F-a_g| \le c_1 ||F-g||_{C^0} \ .
\end{equation}
Now, let us 
consider  the cases  where 
(i) $||F-g|| \ge S/(2c_1)$ 
and (ii) $||F-g|| \le S/(2c_1)$.
In case (i), we get
$$
||RF-Rg||_{C^0}
  \le 2 \le  4 c_1 S^{-1} ||F-g||_{C^0} \ .$$ 
In case (ii), using that  $|a_F| \ge S$  and  $(\ref{qwe233})$, we get
$a_g \ge a_F-S/2 \ge S/2$, and thus, by $(\ref{qwe2})$,
we obtain
$$
\left| a_F^{-1}-a_g^{-1} \right| \le 
a_F^{-1} a_g^{-1}|a_F-a_g| \le  2 S^{-2} c_1 ||F-g||_{C^0} \ .
$$
Hence, again by $(\ref{qwe2})$ and $(\ref{qwe233})$,
there is 
$c_2(N)>0$ with the property that
\begin{eqnarray*}
||RF-Rg||_{C^0}
  &\le &
  ||F^m||_{C^0}|a_F-a_g| +
 |a_g| ||F^m||_{C^1} \left| a_F^{-1}-a_g^{-1} \right| \\
 & & +
  |a_g| ||F^m-g^m||_{C^0} \\
    &\le &  c_2 ||F-g||_{C^0} \ .
  \end{eqnarray*}
Therefore,   this lemma  is satisfied with
$L(N) = \max \{ 2 \epsilon^{-1}, 4 c_1 S^{-1} , c_2 \}$.
\end{proof}

\begin{lemma}
\label{l3}
For all positive constants  $\lambda  < 1$  and $C$
there exist positive constants 
$\alpha(N,\lambda)$ and $n_4(b,N,\lambda,C)$
with the following property:
for every  $f \in \mathcal{I}(N,b)$, and every  $n > n_4$,
 if   $f_n$ is a $C^2$ unimodal map  such that
$$
||R^nf-f_n||_{C^0}< C\lambda^n  \ ,
$$
 then 
  $f_n$ is $[\alpha n+1]$ times renormalizable with $\sigma_{R^mf_n}= \sigma_{R^{n+m}f}$
for every $m=0,\ldots,[\alpha n]$  
(where $[y]$ means the integer part of $y >0$.)
\end{lemma}

\begin{proof}
Let  $\epsilon(N)$  and $n_2(b)$ be as defined in  Lemma  \ref{l2}, 
and let  $L(N)$ and $n_3(b)$  be   as defined in  Lemma   \ref{l1}.
Take $\alpha > 0$     such that
 $L^\alpha \lambda <  1$. 
Set $n_4\ge \max \{n_2, n_3\}$  such that
$C\lambda^{n_4} < \epsilon$ and $C\lambda^{n_4} L^{[\alpha n_4]} < \epsilon$.
Then, for every $n >n_4$, 
  the values 
$C \lambda^n, C \lambda^n L, \ldots,C \lambda^n L^{[\alpha n]}$
are less than $\epsilon$.

 By Lemma \ref{l2}, if $||R^nf-f_n||_{C^0}<   C \lambda^n < \epsilon$
with $n > n_4$, 
then  the 
 map $f_n$ is contained in $D(\sigma_{R^nf})$. 
Thus, $f_n$ is once  renormalizable, and  $\sigma_{f_n}=\sigma_{R^nf}$.
 By induction on  $m=1,\ldots,[\alpha n]$, let us suppose 
 that $f_n$ is $m$ times renormalizable, 
and  $\sigma_{R^i f_n}=\sigma_{R^{n+i}f}$ for every $i=0,\ldots,m-1$. 
By Lemma \ref{l1}, we get that
 $||R^{n+m}f-R^m f_n||_{C^0}< C L^m \lambda^n < \epsilon$. Hence, again by Lemma \ref{l2},  the
 map $R^m f_n$   is once renormalizable, and 
$\sigma_{R^m f_n}=\sigma_{R^{n+m}f}$.
\end{proof}

\begin{lemma}
\label{c3}
There exist positive constants 
$\gamma(N) < 1$, $\alpha(N)$,   $\mu(N)$, and
$C(b,N)$   
with the following property:
for every  $f \in \mathcal{I}(N,b)$,
   there exists    $f_n \in \mathcal{Q(\mu)}$
such that
\begin{rlist}
\item
$||R^nf-f_n||_{C^0} \le C \gamma^n$;
\item  $f_n$ is  $[\alpha n+1]$ times renormalizable  
with  $\sigma_{R^mf_n}= \sigma_{R^{n+m}f}$
for every $m=0,\ldots,[\alpha n]$.
\end{rlist}
\end{lemma}

\begin{proof}
The proof follows from Theorem \ref{th3} and Lemma \ref{l3}.
\end{proof}

\section{Varying quadratic-like maps}

We start by  introducing some classical results on 
Beltrami differentials and holomorphic motions, all of  
which we will apply later in this section 
to vary the combinatorics of quadratic-like maps.

\subsection{Beltrami differentials}

A homeomorphism $h:U \to V$, 
where $U$ and $V$ are contained in $\complexes$ or $\overline{\complexes}$,
 is {\it quasiconformal} 
if it has locally integrable distributional derivatives 
$\partial h$, $\overline{\partial} h$, and if
there is $\epsilon < 1$ with the property that
$\left|\overline{\partial} h/\partial h \right|\le \epsilon$ almost everywhere. 
The Beltrami differential $\mu_h$ of $h$ is
given by 
$\mu_h =\overline{\partial} h/\partial h$.
A quasiconformal   map $h$ is {\it $K$ quasiconformal}   if  
$K \ge (1+||\mu_h||_\infty)/(1-||\mu_h||_\infty)$.

We denote by $D_R(c_0)$  the open disk in $\complexes$
centered at the point $c_0$  and with radius $R>0$.
We also use the notation  $D_R=D_R(0)$ for the disk  centered at the origin.

The following theorem is a slight extension of 
Theorem  4.3 in page 27 of the book  \cite{lehto} by  Lehto.

\begin{theorem}
\label{lehto1}
Let $\psi:\complexes \to \complexes$ be a quasiconformal 
map with the following properties:
\begin{rlist}
\item  $\mu_\psi= \overline{\partial} \psi /\partial \psi$ has  support
contained in the disk $D_R$;
\item $||\mu_\psi||_{\infty} < \epsilon < 1$;
\item $\lim_{|z| \to \infty} (\psi(z)-z)=0$.
\end{rlist}
Then  there exists $C(\epsilon,R)>0$
such that 
$$
||\psi - id||_{C^0} \le C ||\mu_\psi||_{\infty} \ . 
$$
\end{theorem}

\begin{proof}
Let us define $\phi_1=\mu_\psi$, and,  
by induction on  $i \ge 1$, we define $\phi_{i+1}=\mu_\psi H \phi_i$,
where  $H \phi_i$ is the Hilbert  transform of $\phi_i$
given by  the Cauchy Principal Value of 
$$ 
\frac{-1}{\pi} \int \int_{\complexes}
\frac{\phi_i(\xi)}{(\xi-z)^2}
du dv  \ .
$$ 
By Theorem  4.3 in  
page 27 of \cite{lehto},  we get 
$\psi(z)=z+\sum_{i=1}^{\infty} T \phi_i(z)$,
 where $T \phi_i(z)$ is given by
$$ 
\frac{-1}{\pi} \int \int_{\complexes}
\frac{\phi_i(\xi)}{\xi-z}
du dv   \ .
$$ 
By the Calder\'on-Zigmund  inequality (see page 27 of \cite{lehto}),
for every $p \ge 1$,
the Hilbert  operator $H:L^p \to L^p$  is bounded, 
and its norm $||H||_p$  varies  continuously with $p$.
An  elementary integration also shows that $||H||_2=1$ 
(see page 157 of \cite{lehtoVIRT}).
 Therefore, given that $||\mu_\psi||_\infty<\epsilon$,  
there is $p_0(\epsilon)> 2$ with the property  that  
\begin{equation}
\label{azzz1}
||H||_{p_0} ||\mu_\psi||_\infty < ||H||_{p_0} \epsilon < 1  \ .
\end{equation}
Since $p_0>2$, it follows from  H\"older's inequality 
(see page 141 of \cite{lehtoVIRT})
that there is a positive constant $c_1(p_0,R)$ 
such  that  
\begin{equation}
\label{azzz2}
||T \phi_i||_{C^0} \le c_1 ||\phi_i||_{p_0} \ .  
\end{equation}
By a simple  computation, we get 
\begin{equation}
\label{azzz3}
||\phi_i||_{p_0}
\le 
(\pi R^2)^\frac{1}{p_0} ||H||_{p_0}^{i-1}  ||\mu_\psi||_\infty^i \  .
\end{equation}
Thus, by inequalities (\ref{azzz1}), (\ref{azzz2}), and (\ref{azzz3}),
there is a positive constant $c_2(\epsilon,R)$ 
with the property that
\begin{eqnarray*}
||\psi - id||_{C^0} 
& \le &
  \sum_{i=1}^\infty ||T \phi_i||_{C^0}  
\le 
\frac{c_1 (\pi R^2)^\frac{1}{p_0}    ||\mu_\psi||_\infty}  
{1- ||H||_{p_0}  ||\mu_\psi||_\infty} \\
& \le &  c_2  ||\mu_\psi||_\infty  \ .
\end{eqnarray*}
\end{proof}

\subsection{Holomorphic motions}

A {\it holomorphic motion of a subset $X$ of the Riemann sphere over
a disk $D_R(c_0)$}
is a family of maps $\psi_c:X \to X_c$
with the following properties:
(i) $\psi_c$ is an injection of $X$ onto a subset $X_c$ 
of the Riemann sphere;
(ii)   
$\psi_{c_0}=id$;
(iii) for every $z \in X$, 
$\psi_c(z)$ varies holomorphically with $c \in D_R(c_0)$.

\begin{theorem}  
\label{douadyprol} 
{\bf (S\l odkowski \cite{Slodkowski1})}
Let $\psi_c:X \to X_c$
be a holomorphic motion over the  disk $D_R(c_0)$.
Then there is  a holomorphic motion 
$\Psi_c:\overline{\complexes} \to \overline{\complexes}$
over the  disk $D_R(c_0)$ such that
\begin{rlist}
\item  $\Psi_c|X=\psi_c$;
\item  $\Psi_c$ is a $K_c$ quasiconformal map with 
$$
K_c = \frac{R+|c-c_0|}{R-|c-c_0|}  \ .
$$
\end{rlist}
\end{theorem}

 See also Douady's survey \cite{douadyprol11}.

\subsection{Varying the  combinatorics}
\label{llfff}

 Let ${\mathcal M}$ be the set of all  
quadratic-like maps with connected Julia set. 
Let  ${\mathcal P}$ be the set of all normalized quadratic maps $P_c:\complexes \to \complexes$
defined  by  $P_c(z)=1-cz^2$, where $c \in \complexes \setminus \{0\}$.
Two quadratic-like maps $f$ and $g$ are 
{\it hybrid conjugate} if there 
is a quasiconformal conjugacy $h$ between $f$ and $g$
with the property that  $  \overline{\partial} h (z)=0$
for almost every $z \in {\mathcal K}(f)$.  
 By  Douady-Hubbard's Theorem 1 in page 296 of  \cite{douady}, 
for every $f \in {\mathcal M}$ there exists a unique 
quadratic map $P_{c(f)}$ which is hybrid conjugated to $f$.
The map $\xi:{\mathcal M} \to {\mathcal P}$  defined by $\xi(f)=P_{c(f)}$ is called the 
{\it straightening}.

Observe that a real quadratic  map  $P_c$ with $c \notin [1,2]$
has trivial dynamics. Therefore, 
we will restrict our study to the set $\mathcal{Q}([1,2],\mu)$
 of all $f \in \mathcal{Q}(\mu)$ satisfying 
 $\xi(f)=P_{c(f)}$ for some $c(f) \in  [1,2]$.

  Let us choose 
 a radius $\Delta$  large enough 
such  that, for every $c \in [1,2]$,      $P_c(z)=1-cz^2$  is a 
  quadratic-like map when restricted to $P_c^{-1}(D_{\Delta})$.

\begin{lemma}
\label{auxfinding}
There exist positive constants $\Omega(\mu)$ and $K(\mu)$ 
with the following property:
for every    $f \in \mathcal{Q}([1,2],\mu)$
there exists a topological disk $V_f \subset D_{\Omega}$
such that $f$ restricted to $f^{-1}(V_f)$ is a quadratic-like map.
Furthermore, there is a
$K$ quasiconformal homeomorphism 
$\Phi_f:\overline{\complexes} \to \overline{\complexes}$ such that 
\begin{rlist}
\item   $\Phi_f|\Phi_f^{-1}(V_f)$ 
is a hybrid conjugacy between $f$ and $P_{c(f)}$;
\item  $\Phi_f(V_f)=D_\Delta$;
\item $\Phi_f$ is holomorphic over 
$\overline{\complexes} \setminus  \overline{V_f}$; 
\item $\Phi_f(\overline{z})=\overline{\Phi_f(z)}$.
\end{rlist}
\end{lemma}

\begin{proof}
The main point in this proof  is to combine 
the hybrid conjugacy between $f$ and $P_{c(f)}$
given by  Douady-Hubbard, with 
 Sullivan's pull-back argument, and  with McMullen's 
rigidity theorem   for  real quadratic maps. 
Using  Sullivan's pull-back argument  and the hybrid conjugacy 
between $f$ and $P_{c(f)}$, we  
construct a $K$ quasiconformal homeomorphism 
$\Phi_f:\overline{\complexes} \to \overline{\complexes}$ 
 which restricts to a conjugacy 
between $f$ and $P_{c(f)}$.
Moreover, $\Phi_f$ 
satisfies properties (ii), (iii) and (iv)
of this lemma, and  the restriction of $\Phi_f$ to the filled in Julia 
set of $f$ extends to a quasi conformal map that is a hybrid conjugacy 
between $f$ and $P_{c(f)}$. By Rickman's glueing lemma (see Lemma 2 in \cite{douady}) it follows that 
$\Phi_f$ also satisfies  property (i) of this lemma.

Now, we  give the details of the proof:
let us consider  the set of all quadratic-like maps 
 $f:W_f \to W_f'$ contained in  $\mathcal{Q}([1,2],\mu)$.
Using the Koebe Distortion Lemma (see page 84 of \cite{a-bers2}), 
 we can slightly shrink $f^{-n}(W_f')$ for some $n \ge 0$ to 
obtain an open set $V_f$
with the following properties:
\begin{rlist}
\item $V_f$ is
 symmetric with respect to the real axis;
\item the  restriction of $f$ to $f^{-1}(V_f)$ is a 
quadratic-like map;
\item the annulus $V_f \setminus \overline{f^{-1}(V_f)}$ has conformal modulus
 between $\mu/2$ and $2\mu$;
\item 
the boundaries  of $V_f \setminus \overline{f^{-1}(V_f)}$ 
are analytic $\gamma(\mu)$ quasi-circles for some $\gamma(\mu)>0$, 
i. e.,  they are images of an Euclidean circle by   
$\gamma(\mu)$ quasiconformal maps defined on $\overline{\complexes}$. 
\end{rlist}
Let $\mathcal{Q}'$ be the set of all quadratic-like maps 
$f:f^{-1}(V_f)  \to  V_f$ contained in
$\mathcal{Q}([1,2],\mu/2) \cup \mathcal{Q}([1,2],\mu)$
for which  $V_f$ satisfies properties (i), $\ldots$, (iv)
of  last paragraph.
Since for every $f \in \mathcal{Q}'$
the boundaries of $V_f \setminus \overline{f^{-1}(V_f)}$ 
are analytic $\gamma(\mu)$ quasi-circles,
 any convergent sequence $f_n \in \mathcal{Q}'$, with limit $g$, 
in the Carath\'eodory topology has the 
property that the sets $V_{f_n}$ converge 
to $V_g$ in the Hausdorff topology
(see Section 4.1 in  pages 75-76 of \cite{mullen2}). 
Therefore, the set $\mathcal{Q}'$ is closed with respect to the
Carath\'eodory topology,  and hence is  compact. 
Furthermore, by compactness of $\mathcal{Q}'$, and using 
the Koebe Distortion Lemma,
there is an   Euclidean disk $D_\Omega$ 
which contains $V_f$ for every $f \in \mathcal{Q}'$.
 
Now, let us construct  
$\Phi_f:\overline{\complexes} \to \overline{\complexes}$ such that
the properties (i), $\ldots$, (iv) of this lemma are satisfied.

Since $V_f$ is symmetric with respect to the real axis,
there is a unique Riemann Mapping 
$\phi: \overline{\complexes}  \setminus \overline{V_f} \to 
\overline{\complexes} \setminus \overline{D_{\Delta}}$
satisfying   $\phi(\overline{z})=\overline{\phi(z)}$,
and such that $\phi(\reals^+) \subset \reals^+$.
Since the boundaries  of $V_f \setminus \overline{f^{-1}(V_f)}$ 
are analytic $\gamma(\mu)$ quasi-circles,
using the Ahlfors-Beurling Theorem (see Theorem 5.2 in page 33 of \cite{lehto})
 the map $\phi$ has a  
$K_1(\mu)$ quasiconformal homeomorphic extension   
$\phi_1:\overline{\complexes} \to
\overline{\complexes}$ 
which    also is symmetric   $\phi_1(\overline{z})=\overline{\phi_1(z)}$.

Let $\phi_2:V_f \setminus {\mathcal K}(f) \to D_{\Delta} \setminus {\mathcal K}(P_{c(f)})$
 be the unique continuous lift  of $\phi_1$
 satisfying $P_{c(f)} \circ \phi_2(z)= \phi_1 \circ f(z)$,
and such that   $\phi_2(\reals^+) \subset \reals^+$. 
Since $\phi_1$
is a  
$K_1(\mu)$ quasiconformal homeomorphism, so is    $\phi_2$.

Using the Ahlfors-Beurling Theorem, we construct 
a $K_2(\mu)$ quasi-con\-for\-mal homeomorphism
$\phi_3:\overline{\complexes} \setminus {\mathcal K}(f) \to 
\overline{\complexes} \setminus {\mathcal K}(P_{c(f)})$ 
  interpolating $\phi_1$ and $\phi_2$ 
with the following properties: 
\begin{rlist}
\item
 $\phi_3(z) =\phi_1(z)$ 
for every $z \in \overline{\complexes} \setminus V_f$;
\item $\phi_3(z) =\phi_2(z)$ for every $z \in  \overline{ f^{-1}(V_f)}
\setminus  {\mathcal K}(f)$;
\item $\phi_3(\overline{z})=\overline{\phi_3(z)}$.
\end{rlist}
Then the map $\phi_3$ conjugates $f$ on  
$\partial f^{-1}(V_f)$ with $P_{c(f)}$ on  $\partial P_{c(f)}^{-1}(D_\Delta)$, and is   holomorphic over 
$\overline{\complexes} \setminus \overline{V_f} \subset
 \overline{\complexes} \setminus  \overline{D_\Omega}$.

By Theorem 1     in  \cite{douady},   
there is a $K_f'$ quasiconformal 
  hybrid conjugacy $\phi_4:V_f' \to V_{c(f)}'$
between $f$ and $P_{c(f)}$,
where $V_f'$ is a  neigbourhood  of ${\mathcal K}(f)$.
Using the Ahlfors-Beurling Theorem, we construct 
a $K_f''$ quasiconformal homeomorphism
$\Phi_0:\overline{\complexes}   \to \overline{\complexes}$ 
  interpolating $\phi_3$ and $\phi_4$ such that
\begin{rlist}
\item
 $\Phi_0(z) =\phi_3(z)$ 
for every $z \in \overline{\complexes} \setminus f^{-1}(V_f)$;
\item $\Phi_0(z) =\phi_4(z)$ for every $z \in {\mathcal K}(f)$;
\item $\Phi_0(\overline{z})=\overline{\Phi_0(z)}$.
\end{rlist}
Then the map $\Phi_0$   conjugates $f$ on ${\mathcal K}(f) \cup \partial f^{-1}(V_f)$ with $P_{c(f)}$ on ${\mathcal K}(P_{c(f)})  \cup \partial  P_{c(f)}^{-1}(D_\Delta)$,
and  satisfies the properties (ii), (iii) and (iv) as stated in this lemma.
Furthermore, 
$\mu_{\Phi_0}(z)=0$ for   every $z \in \overline{\complexes} \setminus V_f$,
$|\mu_{\Phi_f}(z)| \le (K_2-1)/(K_2+1)$
for a. e. $z \in V_f \setminus f^{-1} ( V_f )$, 
and 
$\mu_{\Phi_f}(z)=0$ for  a. e.  $z \in {\mathcal K}(f) \setminus {\mathcal J}(f)$.

For every $n > 0$,
let us   inductively define
the
$K_f''$ quasiconformal homeomorphism  
$\Phi_n: \overline{\complexes} \to \overline{\complexes}$  as follows:
\begin{rlist}
\item 
$\Phi_n(z)=\Phi_{n-1}(z)$  for every
$z \in \left(\overline{\complexes} \setminus f^{-n}(V_f)\right) \cup {\mathcal K}(f)$;
\item
 $P_{c(f)} \circ  \Phi_n (z)=\Phi_{n-1} \circ f(z)$
 for every $z \in f^{-n}(V_f) \setminus {\mathcal K}(f)$.
 \end{rlist} 
By compactness of the set of all $K_f''$ quasiconformal homeomorphisms on
$\overline{\complexes}$  fixing three points 
 ($0$, $1$ and $\infty$),
 there is a subsequence $\Phi_{n_j}$
which converges to a $K_f''$ quasiconformal homeomorphism
$\Phi_f$. Then $\Phi_f$ satisfies the properties (ii), (iii) and (iv) as stated in this lemma.

The  restriction of 
$\Phi_f$ to the set $f^{-1}(V_f)$ has the property of being a 
quasiconformal conjugacy
between $f$ and $P_{c(f)}$. Furthermore, the Beltrami differential  $\mu_{\Phi_f}$
has the following properties:
\begin{rlist}
\item   
$\mu_{\Phi_f}(z)=0$ for   every $z \in \overline{\complexes} \setminus V_f$;
\item 
$|\mu_{\Phi_f}(z)| \le (K_2-1)/(K_2+1)$
for a. e. $z \in V_f \setminus {\mathcal K}(f)$; 
\item 
$\mu_{\Phi_f}(z)=0$ for  a. e.  $z \in {\mathcal K}(f) \setminus {\mathcal J}(f)$.
\end{rlist}

Therefore, by Rickman's glueing lemma, 
  $\Phi_f:\overline{\complexes} \to \overline{\complexes}$ 
is a $K_2(\mu)$ quasiconformal homeomorphism,
and    $\Phi_f$
restricted to the set $f^{-1}(V_f)$  is a hybrid conjugacy 
between $f$ and $P_{c(f)}$.
\end{proof}

 \bigskip

The lemma below could be proven using the external fibers and the fact that 
the holonomy of the hybrid foliation is quasi conformal as in \cite {lyubich}. 
However we will give a more direct proof of it below.

\begin{lemma}
\label{finding}
There exist positive constants $\beta(\mu) \le 1$, $D(\mu)$, and $\mu'(\mu)$  
with the following property:
for   every  $c \in [1,2]$,
and for every  $f \in \mathcal{Q}([1,2],\mu)$, 
there is  $f_c \in \mathcal{Q}([1,2],\mu')$ 
satisfying $\xi(f_c)=P_c$, and  such that
\begin{equation}
\label{qdq}
||f-f_c||_{C^0(I)}   
\le D |c(f)-c|^\beta  \  .
\end{equation} 
\end{lemma}

\begin{proof}
The main step  of this proof consists of 
  constructing   the real   quadratic-like maps  
$f_c=\psi_c \circ P_c \circ \psi_c^{-1}$ 
satisfying  $f_{c(f)}=f$, and such that the maps
$\omega_c:\overline{\complexes}\to \overline{\complexes}$ 
defined by $\omega_c=\psi_c   \circ \psi_{c(f)}^{-1}$
  form a holomorphic motion    $\omega_c$, and have the property 
of being   holomorphic   on the complement of a disk centered at the origin.
Using Theorem \ref{lehto1}   and  Theorem \ref{douadyprol}, 
we prove that there is  a positive constant $L_3$ 
with the property that $||\omega_c-id ||_{C^0} \le L_3|c-c(f)|$.
  Finally,  we show that this implies the inequality (\ref{qdq}) above.

Now, we   give the details of the proof:
let us choose a small  $\epsilon>0$,  and a small  open set 
   $U$ of $\complexes$ containing the interval $[1,2]$
 such that, for every $c \in U$,  
  the quadratic
map  $P_c(z)=1-cz^2$  has a
  quadratic-like   restriction  to $P_c^{-1}(D_{\Delta})$, and
 $P_c^{-1}(D_{\Delta}) \subset D_{\Delta-\epsilon}$.
Let $\eta: \complexes \to \reals$
 be a $C^\infty$ function with the following properties:
\begin{rlist}
\item $\eta (z)=1$ for every $z \in \complexes \setminus D_{\Delta}$; 
\item $\eta(z)=0$ for every $z \in D_{\Delta-\epsilon}$; 
\item $\eta(z)=\eta(\overline{z})$ for every 
$z \in \complexes$.
\end{rlist}

There is a unique continuous lift 
 $\alpha_{c}:\complexes \setminus P_{c_0}^{-1}(D_{\Delta})
\to \complexes \setminus P_c^{-1}(D_{\Delta})$
of the identity map such that
\begin{rlist}
\item $P_c \circ \alpha_{c}(z)=   P_{c_0}(z)$;
\item $\alpha_{c_0}=id$; 
\item   $\alpha_{c}(z)$ varies continuously with $c$.
\end{rlist}
Then  the  maps
  $\alpha_c$ are  holomorphic injections,  
  and,   for every $z \in \complexes \setminus P_{c_0}^{-1}(D_{\Delta})$,
$\alpha_{c}(z)$ varies holomorphically with $c$.

 Let $\beta_c:\complexes \setminus P_{c_0}^{-1}(D_{\Delta})
\to \complexes \setminus P_c^{-1}(D_{\Delta})$ be 
the   interpolation between the identity map and 
$\alpha_c$ defined by  
 $ \beta_c= \eta \cdot  id + (1-\eta) \cdot \alpha_c$. 
We choose  $r'>0$ small enough such  that,
for every  
$c_0 \in [1,2]$,
and $c \in D_{r'}(c_0) \subset U$,
$\beta_c$ is a diffeomorphism.
Then   
$\beta_c:\complexes \setminus P_{c_0}^{-1}(D_{\Delta})
\to \complexes \setminus P_c^{-1}(D_{\Delta})$
is a holomorphic motion over 
$D_{r}(c_0)$ with the following properties:
\begin{rlist}
\item
the map  $\beta_c$ is a conjugacy between 
$P_{c_0}$ on   
 $\partial P_{c_0}^{-1}(D_{\Delta})$ 
and $P_c$ on  $\partial  P_{c}^{-1}(D_{\Delta})$; 
\item 
the restriction of  $\beta_c$ to the set
$\complexes \setminus D_\Delta$ is the identity map;
\item
if $c$ is real then 
$\beta_c(\overline{z})=\overline{\beta_c(z)}$.
\end{rlist}
 By Theorem \ref{douadyprol}, 
   $\beta_c$ extends to  
a holomorphic motion 
$\hat{\beta}_c:\overline{\complexes} \to \overline{\complexes}$
over    $D_{r'}(c_0)$, and, by taking $r=r'/2$, 
the map    $\hat{\beta}_c$ is $3$ quasiconformal
for every $c \in D_{r}(c_0)$.

By Lemma \ref{auxfinding},  
  there is a $K(\mu)$ quasiconformal homeomorphism 
$\Phi_f:\overline{\complexes} \to \overline{\complexes}$,
and an open set $V_f=\Phi_f^{-1}(D_{\Delta})$ 
such that 
(i)  $\Phi_f$ restricted to
$f^{-1}(V_f)$ 
is a hybrid conjugacy between $f$ and $P_{c(f)}$;
(ii) $\Phi_f$ is holomorphic over $\overline{\complexes} \setminus \overline{V_f}$; and
(iii) $\Phi_f(\overline{z})=\overline{\Phi_f(z)}$.
Let   $\Phi_c:\overline{\complexes} \to \overline{\complexes}$ be defined by
$\Phi_c=\hat{\beta}_c \circ \Phi_f$.
Then,  for every $c \in  D_{r}(c_0)$,  
  $\Phi_c$ is  a $3K$ quasiconformal homeomorphism 
  which conjugates $f$ on  $\partial f^{-1}(V_f)$ 
  with $P_c$ on $\partial P_c^{-1}(D_{\Delta})$.

We define the Beltrami differential $\mu_c$ as follows:
\begin{rlist}
\item $\mu_c(z)=0$ if $z \in {\mathcal K}(P_c) \cup     (\complexes \setminus D_{\Delta})$;
\item $(\Phi_c)^* \mu_c(z)=  0$ if 
$z \in D_{\Delta} \setminus \overline{P_c^{-1} (D_{\Delta})}$; 
\item $ \left(P^n_c\right)_* \mu_c(z)= \mu_c (P^n_c(z))$
if $z \in P_c^{-n} (D_{\Delta}) \setminus \overline{P_c^{-(n+1)} (D_{\Delta})}$ and $n \ge 1$.
\end{rlist}
Then (i) the Beltrami differential 
$\mu_c$  varies holomorphically with $c$;
(ii) $||\mu_c||_{\infty}< (3K-1)/(3K+1)$ for every $c \in D_{r}(c(f))$;
and (iii) if  $c$ is real  then  
$\mu_c(\overline{z})=\overline{\mu_c(z)}$
for almost every $z \in \complexes$.

By the Ahlfors-Bers Theorem (see \cite{a-bers3}),
for every $c \in  D_{r}(c(f))$ there is
a normalized $3K$ quasiconformal homeomorphism   $\psi_c: \overline{\complexes} \to \overline{\complexes}$
with $\psi_c(0)=0$, $\psi_c(1)=1$,  and $\psi_c(\infty)=\infty$ such that $\mu_{\psi_c}=\mu_c$, 
and  $\psi_c(z)$  varies holomorphically
with  $c$.
 Thus,  the restriction of  $\psi_c$     to $\overline{\complexes} \setminus \overline{D_{\Delta}}$ 
is a holomorphic map,  and if $c$ is real then    
$\psi_c(\overline{z})=\overline{\psi_c(z)}$  
for every $z \in \complexes$.

The map $f_c:\psi_c(P_c^{-1}(D_{\Delta})) \to \psi_c(D_{\Delta})$ defined by 
$f_c=\psi_c \circ P_c \circ \psi_c^{-1}$ 
is   $1$ quasiconformal, and thus a holomorphic map. 
Furthermore, the map $f_c$ is hybrid conjugated to $P_c$, and so 
$f_c$ is a  quadratic-like map whose straightening $\xi(f)$ is $P_c$.
Since the conformal modulus of the annulus 
$\psi_c(D_{\Delta}) \setminus \overline{\psi_c(P_c^{-1}(D_{\Delta}))}$
depends only on $3K(\mu)$, we obtain that there is a positive constant $\mu'(\mu)$ 
such that the conformal modulus  of $f_c$ 
is greater than or equal to $\mu'(\mu)$.
If  $c$ is real then  $f_c(\overline{z})=\overline{f_c(z)}$, 
 which implies that $f_c$ is a real quadratic-like map.

For the parameter $c(f)$,   the map $\psi_{c(f)} \circ \Phi_f$ is $1$ quasiconformal
and fixes three points ($0$, $1$ and $\infty$).
Therefore, $\psi_{c(f)} \circ \Phi_f$ is the identity map, and
since the map $\psi_{c(f)} \circ \Phi_f$   conjugates $f$ with $f_{c(f)}$, 
 we get     $f_{c(f)}=f$.

Now, let us prove that the quadratic-like map $f_c$ satisfies inequality
(\ref{qdq}).
By compactness of the set of all 
$3K(\mu)$ quasiconformal homeomorphisms $\phi$ on $\overline{\complexes}$ 
  fixing three  
points ($0$, $1$ and $\infty$), 
    there are positive constants $l(s,\mu) \le s \le L(s,\mu)$    for every $s>0$
with the property that 
\begin{equation}
\label{prr1}
D_l \subset    \phi(D_s)
~~~~~{\rm and}~~~~~
\overline{\complexes} \setminus  \overline{D_{L}} \subset    
\phi(\overline{\complexes} \setminus \overline{D_s})  \  .
\end{equation}
 Thus, there is   ${\Delta}''=L(L({\Delta}))$ 
with the property that 
$\omega_c = \psi_c \circ \psi_{c(f)}^{-1}$ 
is holomorphic in $\overline{\complexes} \setminus \overline{D_{{\Delta}''}}$
 for every $c \in D_{r}(c(f))$,  and $c(f) \in [1,2]$.

Let  $S_{2{\Delta}''}$ be the circle centered at the origin
 and with radius $2{\Delta}''$.
By  (\ref{prr1}), 
we obtain that     $\omega_c(S_{2{\Delta}''})$  
is at a uniform distance from $0$ and $\infty$
for every $c \in D_{r}(c(f))$, and $c(f) \in [1,2]$. 
 Hence,  by the Cauchy Integral Formula,   
and since $\omega_c$ is a holomorphic motion 
over $D_{r}(c(f))$,
the value  $a_c=  \omega_c '(\infty)$  
varies holomorphically 
with $c$,  and there is a constant $L_1(\mu)>0$ 
with the property that
\begin{equation}
\label{333444555}
|a_c-1| < L_1 |c-c(f)| \  .
\end{equation} 
Thus,   (i) the map    $a_c  \omega_c$
is holomorphic in 
$\overline{\complexes} \setminus \overline{D_{\Delta''}}$;
(ii)  $||\mu_{a_c  \omega_c}||_\infty$  
is less than or equal to 
$(9K^2-1)/(9K^2+1)$; and 
(iii) $\lim_{|z| \to \infty} (a_c  \omega_c(z)-z)=0$.
Hence, by Theorem \ref{lehto1},
there is a positive  constant $L_2(\mu)$ such that, 
 for every $c \in D_{r}(c(f))$, and for every $c(f)  \in [1,2]$,
we get
\begin{equation}
\label{333444556}
||a_c  \omega_c-id||_{C^0} \le L_2 ||\mu_{a_c  \omega_c}||_\infty  \  .
\end{equation}
Since $a_c  \omega_c$ is a holomorphic motion 
over $D_{r}(c(f))$, and  by Theorem \ref{douadyprol}, we get
\begin{equation}
\label{pre}
||\mu_{a_c \omega_c} ||_\infty \le \frac{|c-c(f)|}{r} \ .
\end{equation} 
 By inequalities   (\ref{333444555}),  (\ref{333444556}), and (\ref{pre}) 
there is  a positive constant $L_3(\mu)$ such that, 
for every      $c(f) \in [1,2]$,
and for  every $c \in (c(f)-r,c(f)+r)$, we obtain
\begin{equation}
\label{aeq1}
||\omega_c -id||_{C^0(I)} < L_3 |c -c(f)| \ .
\end{equation}
This  implies that 
\begin{equation}
\label{aeq2}
||\omega_c^{-1} -id||_{C^0(I)} < L_3 |c-c(f)| \ .
\end{equation} 
Since $\omega_c$ is a $9K^2$ quasiconformal 
homeomorphism,  and fixes three points,  we obtain from 
Theorem 4.3 in  page 70  of \cite{lehtoVIRT} that  there 
are positive constants $\beta(\mu) \le 1$  and $L_4(\mu)$ 
with the property that
$||\omega_c||_{C^\beta(I)} < L_4$.
Then  by inequalities (\ref{aeq1}) and (\ref{aeq2}) there is 
a positive constant $L_5(\mu)$ such that, 
for every   $c(f) \in [1,2]$,  and for  every $c \in (c(f)-r,c(f)+r)$,
we have 
\begin{eqnarray*}
  ||f_c-f_{c(f)}||_{C^0(I)} 
  & \le  &  ||\omega_c -id||_{C^0(I)}  
  +  
||\omega_c||_{C^\beta(I)}
||P_c-P_{c(f)}||_{C^0(I)} ^\beta \\
&& +
||\omega_c||_{C^\beta}
||P_{c(f)}||_{C^1(I)}^\beta
||\omega_c^{-1} -id||_{C^0(I)}^\beta \\
& \le & L_5|c-{c(f)}|^\beta  \ .
\end{eqnarray*}
Finally, by increasing the constant $L_5$  if necessary,
we obtain that the last inequality is also satisfied
for every   $c(f)$ and $c$ contained  in $[1,2]$. 
\end{proof}

\section{Proofs of the main results}

\subsection{Proof of Lemma \ref{xzx}} 
Let $f=\phi_f \circ p$ be a $C^2$ infinitely renormalizable 
map with bounded combinatorial type.
Let $N$ be such that  
the combinatorial type of $f$ is bounded by $N$,
and set $b=||\phi_f||_{C^2}$.
By Lemma \ref{c3}, 
there are positive constants 
$\gamma(N)<1$, $\alpha(N)$, $\mu(N)$, and $c_1(b,N)$
with the following properties:
  for every $n \ge 0$,
there is an
$[\alpha n+1]$ times renormalizable
quadratic-like map 
$F_n$  with  renormalization type 
$\underline{\sigma}(n)= \sigma_{R^n f},\ldots,\sigma_{R^{n+[\alpha n]} f}$, 
 with conformal modulus greater than or equal to $\mu$,
and satisfying 
\begin{equation}
\label{wwwxxxwww}
||R^nf-F_n||_{C^0(I)} \le c_1\gamma^n \ .
\end{equation}

By Milnor-Thurston's topological classification (see \cite{milnor}
and Theorem 4.2a. in page 470 of \cite{vSdM}),
the real values $c$ 
for which the real quadratic maps $P_c(z)=1-c z^2$ have 
renormalization type $\underline{\sigma}(n)$ is an interval
  $I_{\underline{\sigma}(n)}$.
Thus, by Sullivan's pull-back argument (see \cite{S1}
and Theorem 4.2b. in page 471 of \cite{vSdM}),
there is a unique $c_n \in I_{\underline{\sigma}(n)}$
such that  $P_{c_n}$ has the same combinatorial type  as 
$R^n(f)$.
By Douady-Hubbard's Theorem 1 in  \cite{douady},
there is a unique  quadratic map  $\xi(F_n)=P_{c(F_n)}$ 
which is hybrid conjugated to $F_n$.
Since $F_n$ has renormalization type $\underline{\sigma}(n)$,
the parameter
 $c(F_n)$    belongs to $I_{\underline{\sigma}(n)}$. 
By Lyubich's Theorem 9.6 in page 79 of \cite{lyubich}, there 
are positive constants 
$\lambda(N) < 1$ and $c_2(N)$ 
such that $|I_{\underline{\sigma}(n)}| \le c_2 \lambda^n$.
Therefore, $|c_n-c(F_n)| \le c_2 \lambda^n$.

By Lemma \ref{finding}, there are positive constants 
 $\beta(\mu) <1$, $D(\mu)$, and $\mu'(\mu)$
with the following properties:
for every $n \ge 0$, there is a  real quadratic-like map
$f_n$ with conformal modulus greater than or equal to $\mu'$,
satisfying $\xi(f_n)=P_{c_n}$, and such that
$$
||f_n-F_n||_{C^0(I)} \le D |c_n-c(F_n)|^\beta 
 \le D c_2^\beta \lambda^{\beta n} \ .
$$
Therefore,  the map $f_n$ has the same combinatorial type  as
$R^n(f)$, and, by inequality \eqref{wwwxxxwww}, for 
$C(b,N)= c_1 + D c_2^\beta$  and $\eta (N) =\max\{\gamma,\lambda^\beta\}$,
we get
$$
||R^nf-f_n||_{C^0(I)}    \le   C\eta^n  \ .  
$$ 
\qed

\subsection{Proof of Theorem  \ref{aaa}}
Let $f=\phi_f \circ p$ and $g=\phi_g \circ p$ be any two $C^2$ infinitely renormalizable unimodal maps 
with the same bounded
  combinatorial type. 
Let $N$ be such that  
the combinatorial type of $f$ and $g$ are bounded by $N$,
and set $b=\max\{||\phi_f||_{C^2}, ||\phi_g||_{C^2}\}$.
For  every $n \ge 0$,
let $m= [\alpha n]$, where $0<\alpha < 1$ 
will be fixed  later in the proof.
 By  Lemma \ref{xzx}, there are positive constants  
  $\eta(N)<1$ and  $c_1(b,N)$, 
 and   there are 
infinitely renormalizable real quadratic-like maps
 $F_m$ and $G_m$  with the following property:  
\begin{equation}
\label{qq55}
||R^mf-F_m||_{C^0(I)} \le
  c_1 \eta^{\alpha n}
~~~~~{\rm and}~~~~~
||R^mg-G_m||_{C^0(I)} \le
  c_1 \eta^{\alpha n} \ .
\end{equation} 
By Lemma \ref{l1},
there are positive constants  $n_3(b)$ and  $L(N)$
such  that, for every $m>n_3$, we get
\begin{eqnarray}
\label{qq11}
||R^{n}f-R^{n-m}F_m||_{C^0(I)}
 & \le &
 L^{n-m}||R^{m}f-F_m||_{C^0(I)} \\
 & \le &
c_1 \left( L^{1-\alpha}  \eta^\alpha \right)^n  \ ,  \nonumber 
\end{eqnarray}
and, similarly,
\begin{equation}
\label{qq11111}
||R^{n}g-R^{n-m}G_m||_{C^0(I)} \le 
c_1 \left( L^{1-\alpha}  \eta^\alpha \right)^n  \  . 
\end{equation}
Now, we  fix   $0 < \alpha(N) < 1$   
such that  $L^{1-\alpha}  \eta^\alpha  <1$.

Again, by  Lemma \ref{xzx}, 
 $F_m$ and $G_m$ have conformal modulus greater than or equal to $\mu(N)$,
and  the same combinatorial type 
 as $R^mf$ and $R^mg$.
Therefore, by McMullen's Theorem 9.22 in page 172 
 of \cite{mullen2},
there are positive constants   $\nu_2(N)<1$ and $c_2(\mu,N)$
with the property that
\begin{equation}
\label{r3}
||R^{n-m}F_m-R^{n-m}G_m||_{C^0(I)} \le c_2 \nu_2^{n-m} \ .
\end{equation}
By inequalities $(\ref{qq11})$,  $(\ref{qq11111})$, 
and  $(\ref{r3})$, there are constants $c_3 (b,N) = 2c_1 + c_2$ and
 $\nu_3 (N) = \max \{L^{1-\alpha}  \eta^\alpha ,\nu_2^{1-\alpha}\}$ such that
$$
||R^{n}f-R^{n}g||_{C^0(I)} \le
 c_3 \nu_3^{n}   \  .
$$
By   Theorem 9.4 in page 552 of 
\cite{vSdM}, 
the exponential convergence implies that
there is a $C^{1+\alpha}$ diffeomorphism
which conjugates  $f$ and $g$
along the closure of the corresponding  orbits of the critical points
 for some $\alpha(N) > 0$.
\qed

\bigskip\bigskip\bigskip 
The exponential convergence of the renormalization operator in the space 
of real analytic unimodal maps holds for every combinatorial type.  
Indeed, if $f$ and $g$ are real analytic infinitely renormalizable 
maps, by the complex bounds 
 in Theorem A of   Levin-van Strien in \cite{levin}, 
there 
exists an integer $N$ such that $R^N(f)$ and $R^N(g)$ have quadratic 
like extensions. Then we can use Lyubich's Theorem 1.1 in 
\cite{lyubich3}  to conclude the exponential convergence. 
However, as we pointed out before, this is not sufficient to give the 
$C^{1+\alpha}$ rigidity. Finally, at the moment, we cannot prove the 
exponential convergence of the operator for $C^2$ mappings with unbounded 
combinatorics.

\par\bigskip
\par\bigskip
\noindent {\large \bf Acknowledgments}
Alberto Adrego Pinto would like to thank 
IMPA, University of Warwick,  and IMS at SUNY Stony Brook
for their hospitality. 
We would like to thank Edson de Faria, and  Mikhail 
Lyubich  for useful discussions. 
This work has been partially supported by 
the Pronex Project on Dynamical Systems, 
 Funda\c{c}\~ao para a Ci\^encia, 
Praxis XXI,   Calouste Gulbenkian Foundation,  and
  Centro de Matem\'atica Aplicada, 
 da
Universidade do Porto, Portugal.

\end{document}